\documentclass{gen-j-l}%
\usepackage{amssymb}
\usepackage{amsmath}
\usepackage{amsfonts}
\usepackage{graphicx}%
\newtheorem{theorem}{Theorem}[section]
\newtheorem{lemma}[theorem]{Lemma}
\theoremstyle{definition}
\newtheorem{definition}[theorem]{Definition}

\theoremstyle{remark}

\numberwithin{equation}{section}
\theoremstyle{plain}

\newtheorem{conjecture}{Conjecture}

\newtheorem{proposition}{Proposition}

\copyrightinfo{2001}{enter name of copyright holder}
\begin{document}
\title[Harmonic Maps and Bubbling at Boundary of Moduli Space]{Harmonic Maps of Surfaces Approaching the Boundary of Moduli Space and
Eliminating Bubbling}
\author{Simon P. Morgan}
\address{University of Minnesota}
\email{morgan@math.umn.edu}
\urladdr{http://www.math.umn.edu/\symbol{126}morgan}
\keywords{Harmonic maps, bubbling, dimension collapsing}

\begin{abstract}
The limit of energies of a sequence of harmonic maps as their annular domains
approach the boundary of moduli space depends upon the boundary point
approached. The infinite energy case is associated with limits of images
containing ruled surfaces. The finite energy case yields a limit of images,
under a suitable topology, with a union of discs and straight line segments.
Generalization to higher numbers of boundary components shows that minimal
surfaces union straight line segments can still be achieved, and that the
configuration of straight line segments depends on the direction of approach
of domain conformal classes to the boundary point of domain moduli space.
Bubbling can occur in variable ways according to the metric representatives of
the conformal classes of domains. A method is given whereby bubbling can be
eliminated yielding a point-wise limit of maps with image limit containing a
surface union straight line segments.

\end{abstract}
\maketitle
\tableofcontents

\section{Introduction}

\subsubsection{\textbf{Example for new compactness theorems}}

This paper, the first of two with [MS1], provides an example (Figure 1) of a
physical optimization problem involving minimal surfaces and threads of
viscoelastic fluids with prescribed boundary and initial conditions.

\begin{center}%
%TCIMACRO{\FRAME{itbpF}{0.8867in}{0.8916in}{0in}{}{}{Figure}%
%maintain-aspect-ratio TRUE;  display "USEDEF";  valid_file "T";
%width 0.8867in;  height 0.8916in;  depth 0in;  original-width 1.0726in;
%original-height 1.0793in;  cropleft "0.007492";  croptop "1.052154";
%cropright "1.007492";  cropbottom "0.052154";
%tempfilename '../Sphere bundle compactness/I3LT5O02.wmf';tempfile-properties "XPR";}%
%}}%
%BeginExpansion
{\includegraphics[
trim=0.008036in 0.056290in -0.008036in -0.056290in,
%natheight=1.079300in,
%natwidth=1.072600in,
%height=0.8916in,
width=0.8867in
]%
{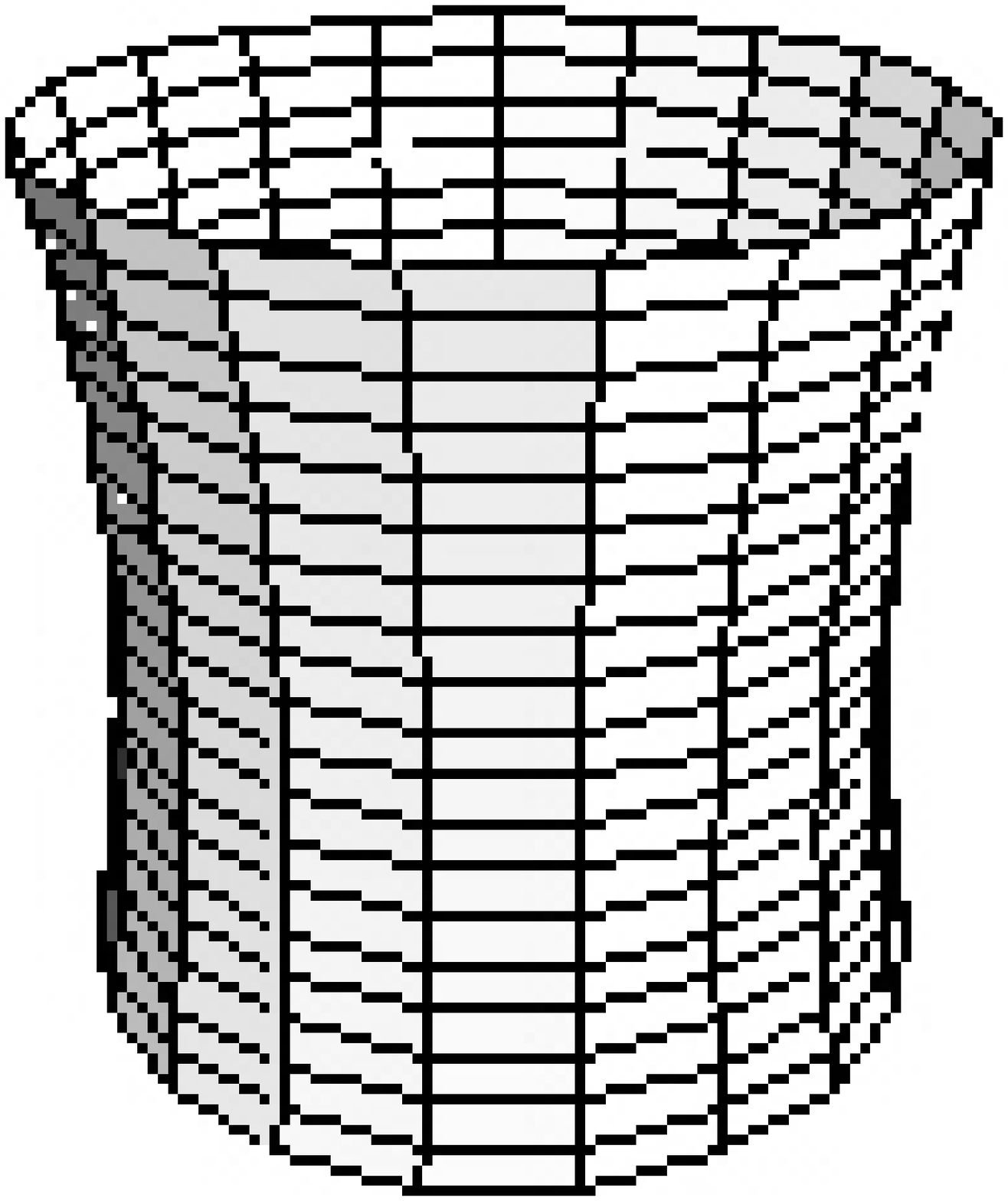}%
}%
%EndExpansion%
%TCIMACRO{\FRAME{itbpF}{0.8867in}{0.895in}{0in}{}{}{Figure}%
%maintain-aspect-ratio TRUE;  display "USEDEF";  valid_file "T";
%width 0.8867in;  height 0.895in;  depth 0in;  original-width 1.0726in;
%original-height 1.0826in;  cropleft "0";  croptop "1";  cropright "1";
%cropbottom "0";  tempfilename 'I3LT5O03.wmf';tempfile-properties "XPR";}}}%
%BeginExpansion
{\includegraphics[
%natheight=1.082600in,
%natwidth=1.072600in,
%height=0.895in,
width=0.8867in
]%
{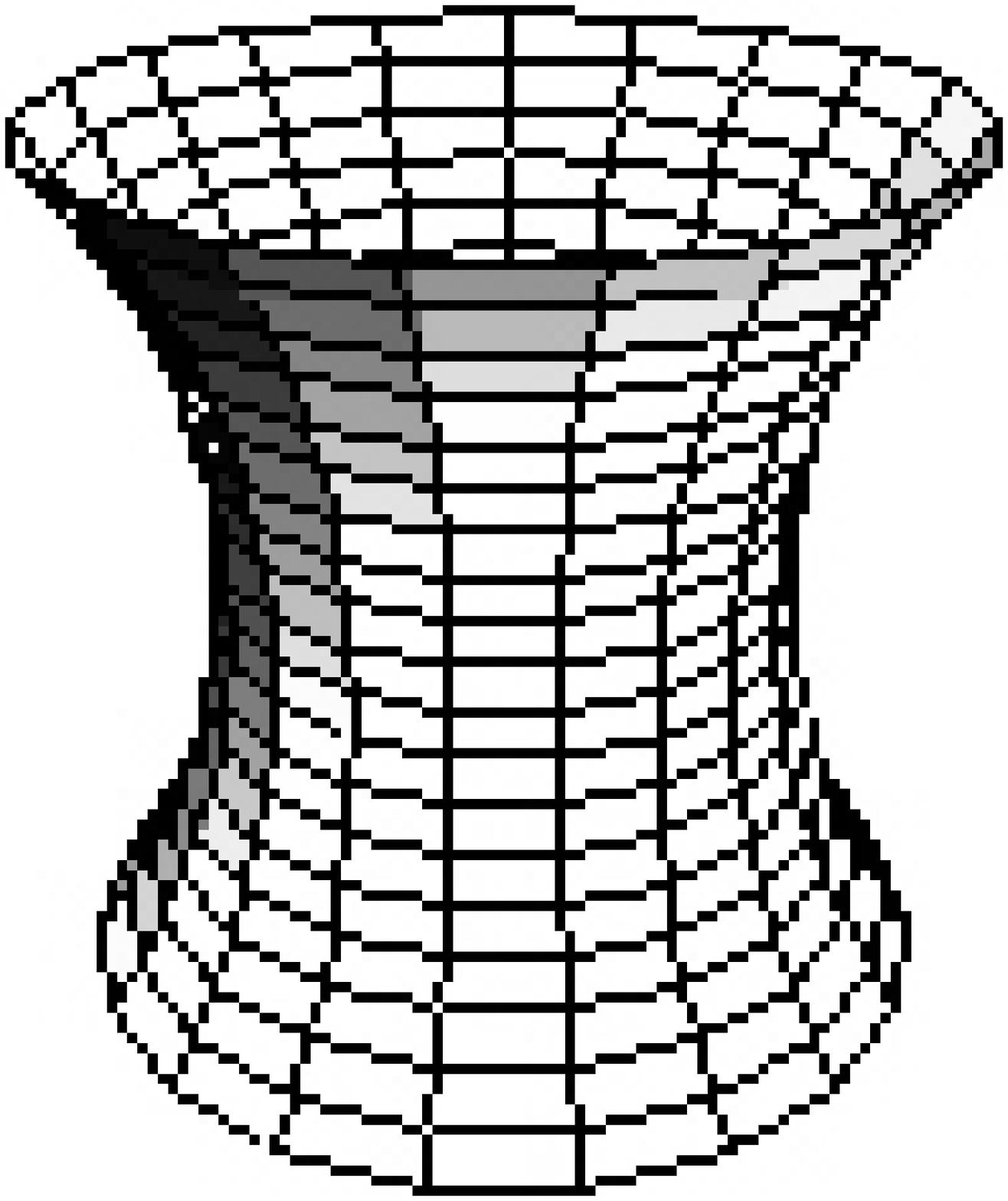}%
}%
%EndExpansion%
%TCIMACRO{\FRAME{itbpF}{0.8867in}{0.8916in}{0in}{}{}{Figure}%
%maintain-aspect-ratio TRUE;  display "USEDEF";  valid_file "T";
%width 0.8867in;  height 0.8916in;  depth 0in;  original-width 1.0726in;
%original-height 1.0793in;  cropleft "0";  croptop "1";  cropright "1";
%cropbottom "0";  tempfilename 'I3LT5O04.wmf';tempfile-properties "XPR";}}}%
%BeginExpansion
{\includegraphics[
%natheight=1.079300in,
%natwidth=1.072600in,
%height=0.8916in,
width=0.8867in
]%
{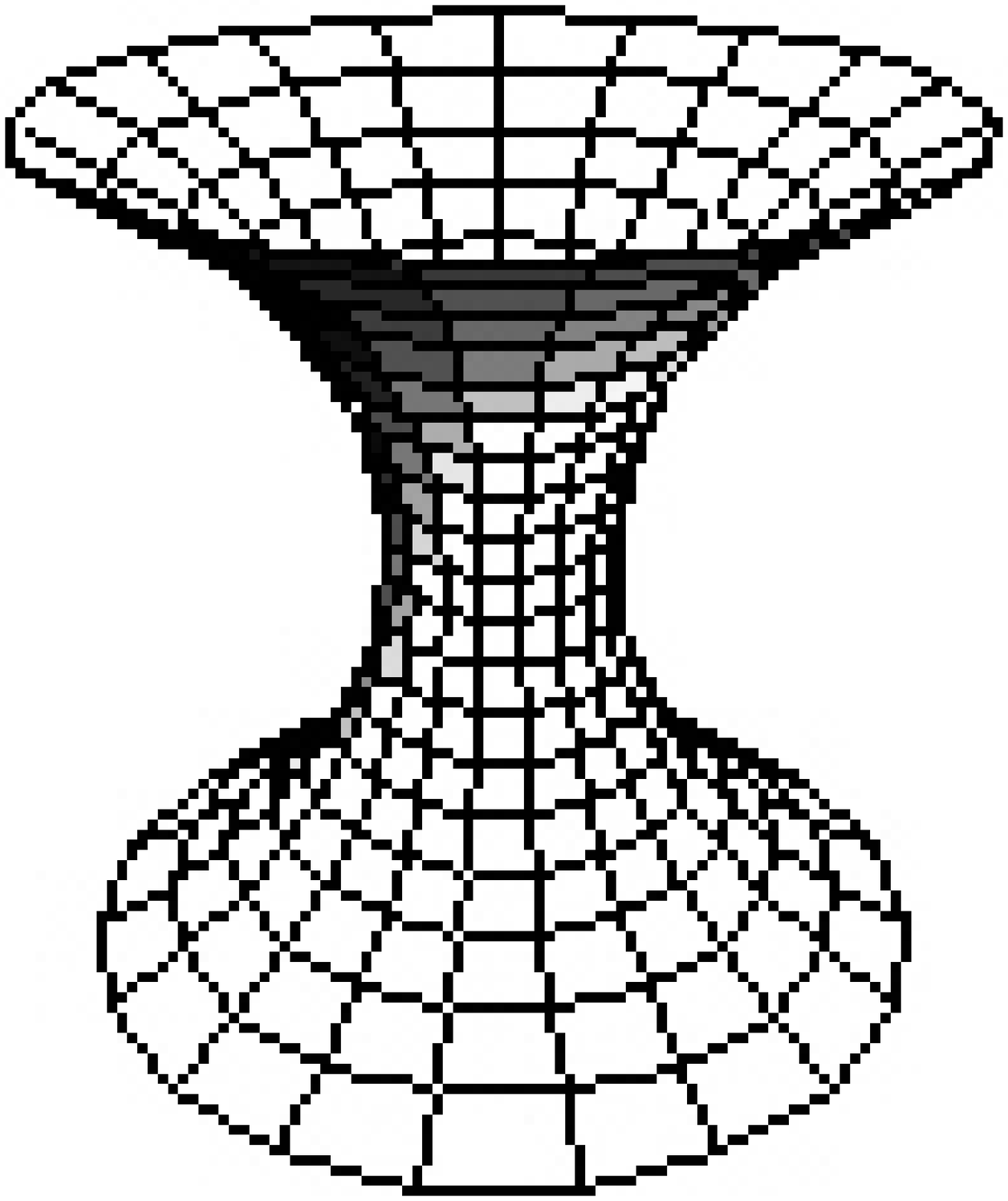}%
}%
%EndExpansion%
%TCIMACRO{\FRAME{itbpF}{0.8867in}{0.8916in}{0in}{}{}{Figure}%
%maintain-aspect-ratio TRUE;  display "USEDEF";  valid_file "T";
%width 0.8867in;  height 0.8916in;  depth 0in;  original-width 1.0726in;
%original-height 1.0793in;  cropleft "0";  croptop "1";  cropright "1";
%tempfilename '../Sphere bundle compactness/I3LT5O05.wmf';tempfile-properties "XPR";}%
%}}%
%BeginExpansion
{\includegraphics[
%natheight=1.079300in,
%natwidth=1.072600in,
%height=0.8916in,
width=0.8867in
]%
{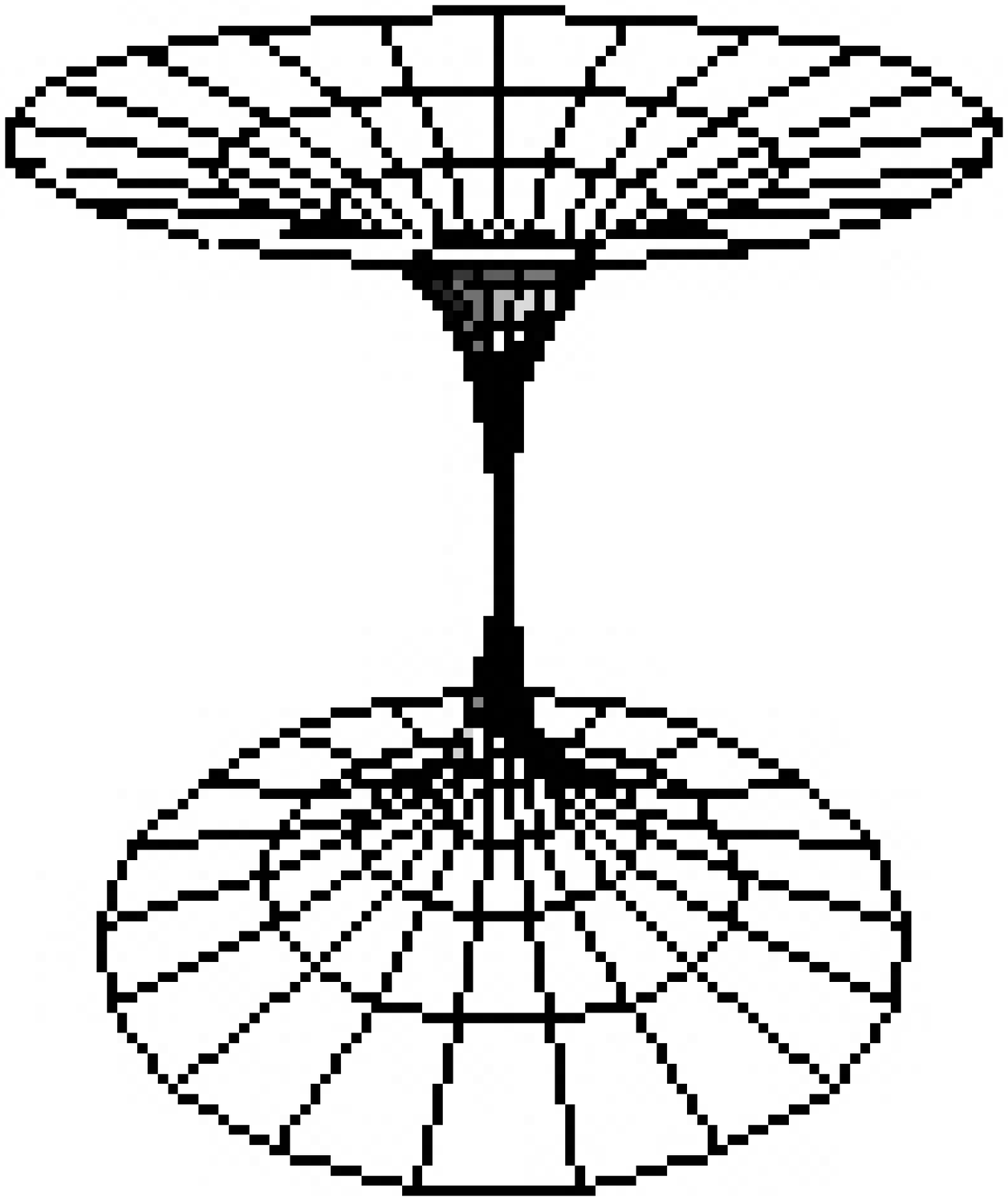}%
}%
%EndExpansion%
%TCIMACRO{\FRAME{itbpF}{0.8925in}{0.8983in}{0in}{}{}{Figure}%
%maintain-aspect-ratio TRUE;  display "USEDEF";  valid_file "T";
%width 0.8925in;  height 0.8983in;  depth 0in;  original-width 1.0809in;
%original-height 1.0867in;  cropleft "0";  croptop "1";  cropright "1";
%tempfilename '../Sphere bundle compactness/I3LT5O06.wmf';tempfile-properties "XPR";}%
%}}%
%BeginExpansion
{\includegraphics[
%natheight=1.086700in,
%natwidth=1.080900in,
%height=0.8983in,
width=0.8925in
]%
{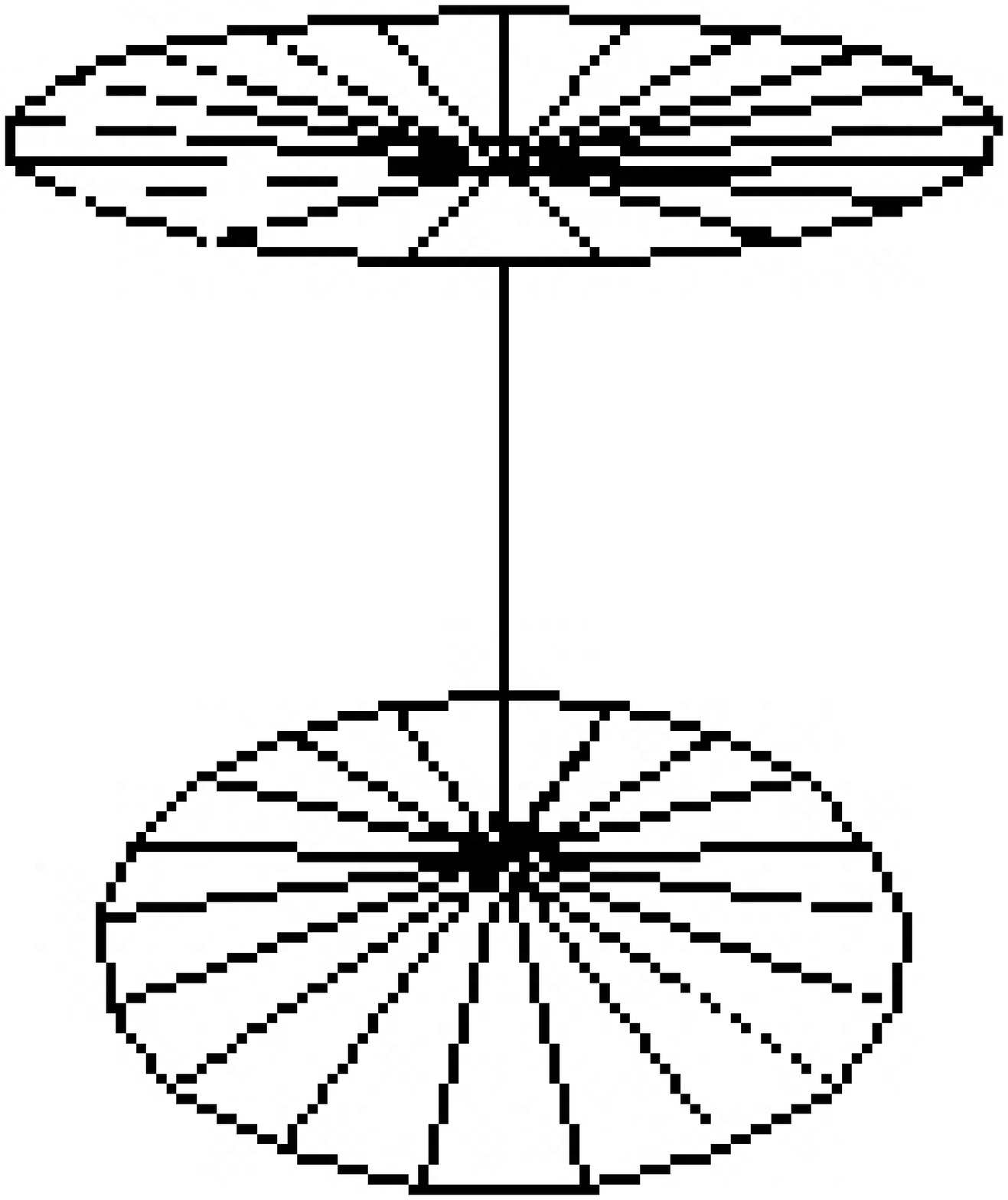}%
}%
%EndExpansion

\textbf{Figure 1: Dimension collapsing}\bigskip
\end{center}

This has a solution which is the limit of minimizers of an energy functional,
that is the limit of images of harmonic maps. The regularity associated with
minimizers, such as harmonic map images, can help establish existence of
limits of such minimizers. The limiting process in our example has two
features of interest. First it involves dimension collapsing of the images of
the harmonic maps where the lower dimensional sets are part of the solution,
under a suitable topology. Secondly discontinuities can arise in bubbling.
Although in our examples, each limit of images can be achieved as a continuous
point-wise limit by careful construction of a sequence of harmonic maps, in
general a sequence will not have a continuous point-wise limit.

The second paper [MS1] provides two versions of a topology that overcome both
the problem of requiring the lower dimensional sets to be kept in the limit,
and the problem of discontinuities due to bubbling in a sequence of harmonic
map images. This is done using sphere bundle measures (in $\mathbb{R}^{n}%
$\textsf{X}$\mathbb{S}^{n-1}$) to represent sets in $\mathbb{R}^{n}$ and take
limits. The sphere bundle measure compactness has the added advantage of
providing a topology which compactifies sufficiently regular unions of
rectifiable subsets of $\mathbb{R}^{n}.$ Fortunately our motivating example
consists of taking limits of images of harmonic maps which have sufficient regularity.

\subsubsection{\textbf{Minimal surfaces union straight line segments}}

Douglas [D] showed that minimal surfaces which are topological discs can be
images of harmonic maps of fixed discs with some variability in the boundary
conditions. We can generalize this method to genus zero minimal surfaces with
more than one boundary component by not only varying the boundary conditions,
but by varying the conformal class of the domain too.

Consider the case of two boundary components. We take energy minimizing
sequences of harmonic maps over varying annular domains. This can result in
minimal surfaces union straight line segments as limits of images (figure 1).
Such collapsing of surfaces down to line segments are mentioned in the context
of harmonic maps [SU] as `bridges' between minimal spheres. In [Fo][DF], where
this example is discussed, extraordinary homology is used to define a homology
class of solutions that can contain not only surfaces but the limit of
surfaces as parts of them collapse down to line segments.

\subsubsection{\textbf{Harmonic map existence.}}

The original motivation for the work reported here was to study harmonic maps
of surfaces which have singular images such as cones. In our example these are
image limits as the other end of moduli space is approached. This can also be
thought of as an inverse problem: What domains can map harmonically onto a
given image embedded in Euclidean space? Harmonic maps are solutions to an
energy minimizing variational problem with fixed range, domain and boundary
conditions. There are many regularity and existence results for harmonic maps,
e.g.: [SCU][EL1][EL2][EF].

Sometimes solutions do not exist, even singular ones. Say we wish to achieve
an embedded cone as an image of a harmonic map from a disc into $\mathbb{R}%
^{3}$. If we place boundary conditions on a flat domain disc in $\mathbb{R}%
^{3}$ to be the identity, and also require the center of the disc to be mapped
to a point out of the plane, any sequence maps whose energy tends to the
infimum will have a point-wise limit with at least one discontinuity as in
Figure 2.%

%TCIMACRO{\FRAME{dtbpFU}{2.8825in}{0.8966in}{0pt}{\Qcb{\textbf{Figure 2:
%Harmonic map non-existence}}}{}{Figure}{\special{ language "Scientific Word";
%type "GRAPHIC";  maintain-aspect-ratio TRUE;  display "USEDEF";
%valid_file "T";  width 2.8825in;  height 0.8966in;  depth 0pt;
%original-width 2.7538in;  original-height 0.8385in;  cropleft "0";
%tempfilename 'I3LT5O07.wmf';tempfile-properties "XPR";}}}%
%BeginExpansion
\begin{center}
\includegraphics[
%natheight=0.838500in,
%natwidth=2.753800in,
%height=0.8966in,
width=2.8825in
]%
{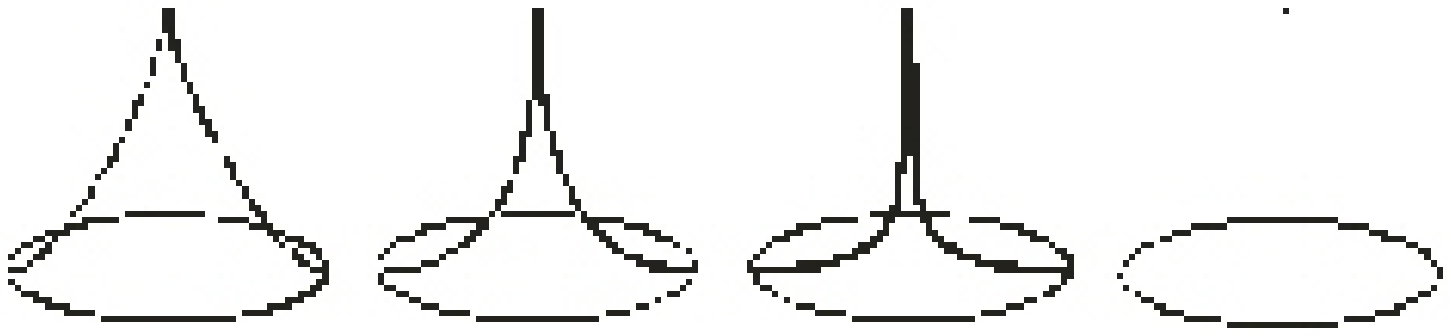}%
\\
\textbf{Figure 2: Harmonic map non-existence}%
\end{center}
%EndExpansion

Overcoming the non-existence of any harmonic maps for this problem motivated
studying maps from annular domains rather than discs. The outer boundary can
be mapped to the disc boundary in the range and the inner boundary can be
mapped to the point out of the plane. Harmonic maps do exist for this modified
problem. However a ruled surface cone can only be achieved as a limit of
images of harmonic maps whose energy tends to infinity. See table 1.

\subsubsection{\textbf{Techniques and bubbling}}

We will use two classical results specific to harmonic maps of surfaces.
Firstly we use the energy-area inequality for maps of surfaces, $Energy$
$\geq$ $2$($Area$ $of$ $Image)$ and an improved version (2.4) using image
curvature [MS2].

Secondly we use the fact that the image of a harmonic map of a surface is an
invariant of the conformal structure of the domain. Any two harmonic maps of
surfaces from conformally equivalent domains with compatible boundary
conditions will have the same image and the same energy. In our annular domain
case, this gives us a one dimensional moduli space of conformal structures
which when compactified has two boundary points. This gives us a corresponding
1-parameter family of images with two limits as shown in table 1. See [Ah] and
[N] for accounts of the moduli space of multiply connected planar domains.

This invariance of image and energy with respect to domain conformal class
means that there are two types of change to domains for harmonic maps of
surfaces. Changes of conformal structure, that will in general change the
image and energy of a map, and then changes of metric structure within a fixed
conformal class. The latter type of changes do not change the image, but can
introduce or eliminate the phenomenon of bubbling in sequences of harmonic maps.

We will define bubbling as the phenomenon of a sequence of maps from a fixed
domain to a range having a point set limit as a graph but may not necessarily
have a point-wise limit as a function because the graph becomes vertical for a
set of positive measure in the range. For example the limit as n$\rightarrow
\infty$ of the graphs of $f(x)=x^{1/n}$ contains the unit interval on the Y
axis. See section 3 for an example of bubbling with point-wise limits. [SU]
gave an early example of bubbling in terms of fractional linear
transformations on spheres being conformal, and hence harmonic, but not
compact. [W] reviews a variety of research areas where bubbling has to be
dealt with.

Consider as an example of bubbling, the sequence of maps from the compactified
complex plane to itself, $z\rightarrow nz$. The image of the point-wise limit
as $n\rightarrow\infty$ is $\{\infty\}$ $\cup\{0\}$, however the Hausdorff set
limit of images is the compactified complex plane. This poses problems for
trying to study limits of images of maps as the image of the point-wise limit
of a map, as parts of the image can be lost in the limit. Also it may lead to
a change in the degree of maps between the limits of images and images of
point-wise limits of maps. We eliminate bubbling in section 3 similarly to
[SU] where parameterizations are chosen to be `in balance'.

The development of minimal surface theory after Douglas included many
mathematical approaches, and we shall use one of these in particular,
geometric measure theory [S][MF][F]. Federer and Fleming [FF] used current
compactness to prove the existence of minimal surfaces in $\mathbb{R}^{n}$
with certain fixed boundary. We will use varifold and current compactness in a
novel way on the sphere bundle [MS1] to take limits of the images. This avoids
the problems of bubbling for point-wise limits, and captures lower dimensional
sets in the limit. The Hausdorff set topology also has these advantages, but
does not ensure the regularity, rectifiability and measurability, given by
geometric measure theory.

\section{The annulus example}

\subsection{The one parameter families of domains, images and their limits}%

%TCIMACRO{\FRAME{dtbpFU}{4.8186in}{3.8539in}{0pt}{\Qcb{\textbf{Table 1}}}%
%{}{Table 1}{\special{ language "Scientific Word";  type "GRAPHIC";
%maintain-aspect-ratio TRUE;  display "USEDEF";  valid_file "T";
%width 4.8186in;  height 3.8539in;  depth 0pt;  original-width 6.0548in;
%original-height 4.6185in;  cropleft "0";  croptop "1.0470";  cropright "1";
%cropbottom "0";  tempfilename 'I3LT5O08.wmf';tempfile-properties "XPR";}}}%
%BeginExpansion
\begin{center}
\includegraphics[
trim=0.000000in 0.000000in 0.000000in -0.217070in,
%natheight=4.618500in,
%natwidth=6.054800in,
%height=3.8539in,
width=4.8186in
]%
{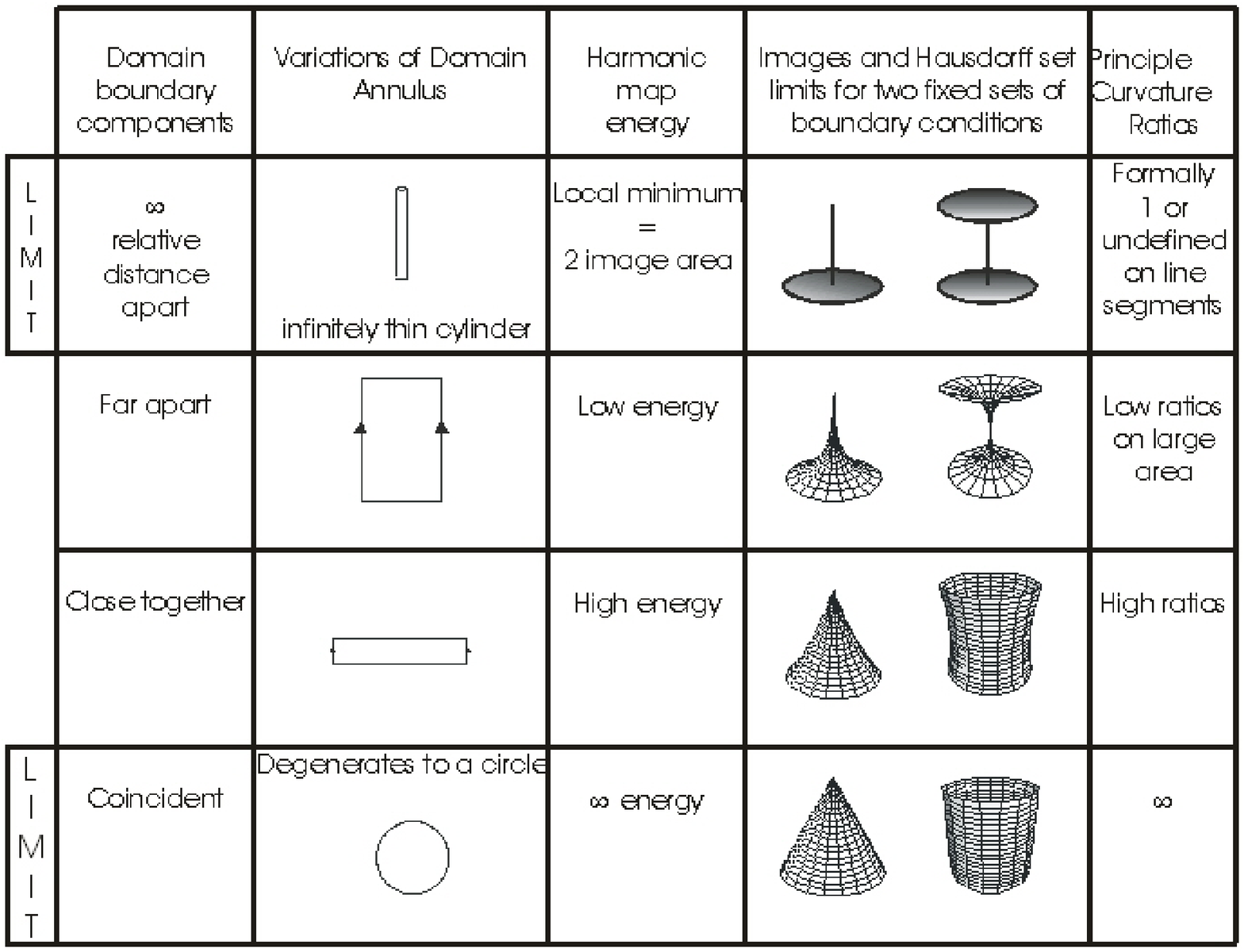}%
\\
\textbf{Table 1}%
\end{center}
%EndExpansion

Table 1 shows two domain annuli for a harmonic map and their respective images
(rows 2 and 3). A geometric interpretation of the two limiting cases of
domains and images is given (rows 1 and 4). Note that the energy-area
inequality for harmonic maps is true in row 1 of the table, where we do not
have a harmonic map but a limit of energies of harmonic maps and a limit of
their images. In row 4, the limit of energies is trivially greater than twice
the areas of the limits of images. Using (2.4) the first inequality has
infinity on both sides for this case.

The images for domains on the interior of moduli space are obtained as
solutions (2.2) to an ODE (2.1) based on the first variation of energy of a
radially symmetric map $h(x,\theta)=(R(x),\Theta\left(  \theta\right)
,Z(x))$, of an annulus based on identified rectangle of width $2\pi/a$,
$(0<\theta<2\pi)$, and height a , ($-a/2<x<a/2$). For each row in table 1 we
can give values of $a$. See table 2.

\begin{center}%
\begin{tabular}
[c]{|c|c|}\hline
Row in Table 1 & \textit{a}\\\hline
1 & $\infty$\\\hline
2 & 3\\\hline
3 & 0.5\\\hline
4 & 0\\\hline
\end{tabular}

\textbf{Table 2}
\end{center}

Note: $a\in\left(  0,\infty\right)  $ is a parameter representing the
conformal class of the domain annulus.

Let $\varsigma\left(  x\right)  $ be any smooth function with compact support
that vanishes on the boundary. We take the derivative of the energy of the
deformed family of maps $h(x,\theta)=(R(x)+\varsigma\left(  x\right)
t,\Theta\left(  \theta\right)  ,Z(x))$, with respect to time and set it equal
to zero. Note we are only considering radially symmetric deformations only.
The energy can be expressed in terms of squares of the entries of the Jacobian
in local orthonormal coordinates yielding:

\begin{center}%
\[
\frac{d}{dt}\mid_{t=0}\int\left(  \left(  R_{x}+\varsigma_{x}(x)t\right)
^{2}+\left(  \frac{\partial\Theta}{\partial\theta}\right)  ^{2}\left(
R+\varsigma(x)t\right)  ^{2}\right)  dA=0
\]

\end{center}

After multiplying out, differentiating under the integral sign, using
integration by parts, and setting the integrand to zero yields:

\begin{center}%
\[
a^{2}R\varsigma(x)-R_{xx}\varsigma(x)=0
\]

This is always satisfied when \
\begin{equation}
a^{2}R-R_{xx}=0
\end{equation}

\end{center}

The general solution for the image in cylindrical $R,\Theta,Z$ coordinates is
of the form;%

\begin{equation}
R=\frac{A\cosh\left(  a^{2}Z\right)  }{\cosh\left(  a^{2}\right)  }%
+\frac{B\sinh\left(  a^{2}Z\right)  }{\sinh\left(  a^{2}\right)  }%
\end{equation}

where $a\in\left(  o,\infty\right)  $ is a parameter representing the
conformal class of the domain annulus. Also similar calculations show that $Z
$ turns out to be a linear function of $x$, and that the radially symmetric
map is energy stationary.

Certain of these hyperbolic catenoids ($B=0$) will be minimal surfaces (see
[MF][Fo][DF] for examples), and these may be stable or unstable. For fixed
boundary conditions such minimal surfaces occur only if the boundary circles
are sufficiently close together for the given radii.

In higher dimensions harmonic maps from products of spheres and an interval
give an analogous solution. A product of an $\mathbb{S}^{n}$ of radius$\frac
{1}{^{n}\sqrt{a}}$ with an interval of length $a$, will map harmonically into
$\mathbb{R}^{n+2}$ with unit $\mathbb{S}^{n}$ at a distance 1 apart as
boundary conditions. Equation (2.2) now generalizes to%

\begin{equation}
R=\frac{A\cosh\left(  na^{1+\frac{1}{n}}Z\right)  }{\cosh\left(
na^{^{1+\frac{1}{n}}}\right)  }+\frac{B\sinh\left(  na^{^{1+\frac{1}{n}}%
}Z\right)  }{\sinh\left(  na^{^{1+\frac{1}{n}}}\right)  }%
\end{equation}

This means that in higher dimensions the same qualitative behavior as depicted
in table 1 occurs. The image becomes a cylinder at one limit and a union of a
straight line segment and two $n+1$ balls at the other limit. Note also that
suitably interpreted, equation 2.3 holds for $n=0$, giving $R=A$.

\subsection{Physical intuition I: Lengths of geodesics and domain conformal
structure}

The moduli space of a domain surface ([Ah][N]), is the space of all
conformally equivalent structures on the surface. In the case of the domain
annuli in table 1, the ratio of height to width of the identified rectangle
indicates the conformal structure. The moduli space of the torus minus two
discs in figure 3 can be considered in terms of the ratio of lengths of curves
connecting distinct pairs of boundary components (a) and curves going round
homologically distinct and non-trivial loops ($b,c,d,e$). Note that some of
these can be the boundary curves ($d$), ($f$).

\begin{center}%
%TCIMACRO{\FRAME{dtbpFU}{2.4201in}{1.5683in}{0pt}{\Qcb{\textbf{Figure 3: Curve
%lengths represent position in moduli space}}}{}{Figure}%
%maintain-aspect-ratio TRUE;  display "USEDEF";  valid_file "T";
%width 2.4201in;  height 1.5683in;  depth 0pt;  original-width 6.7878in;
%original-height 4.3694in;  cropleft "0";  croptop "1";  cropright "1";
%cropbottom "0";  tempfilename 'I3LT5P09.wmf';tempfile-properties "XPR";}}}%
%BeginExpansion
\begin{center}
\includegraphics[
%natheight=4.369400in,
%natwidth=6.787800in,
%height=1.5683in,
width=2.4201in
]%
{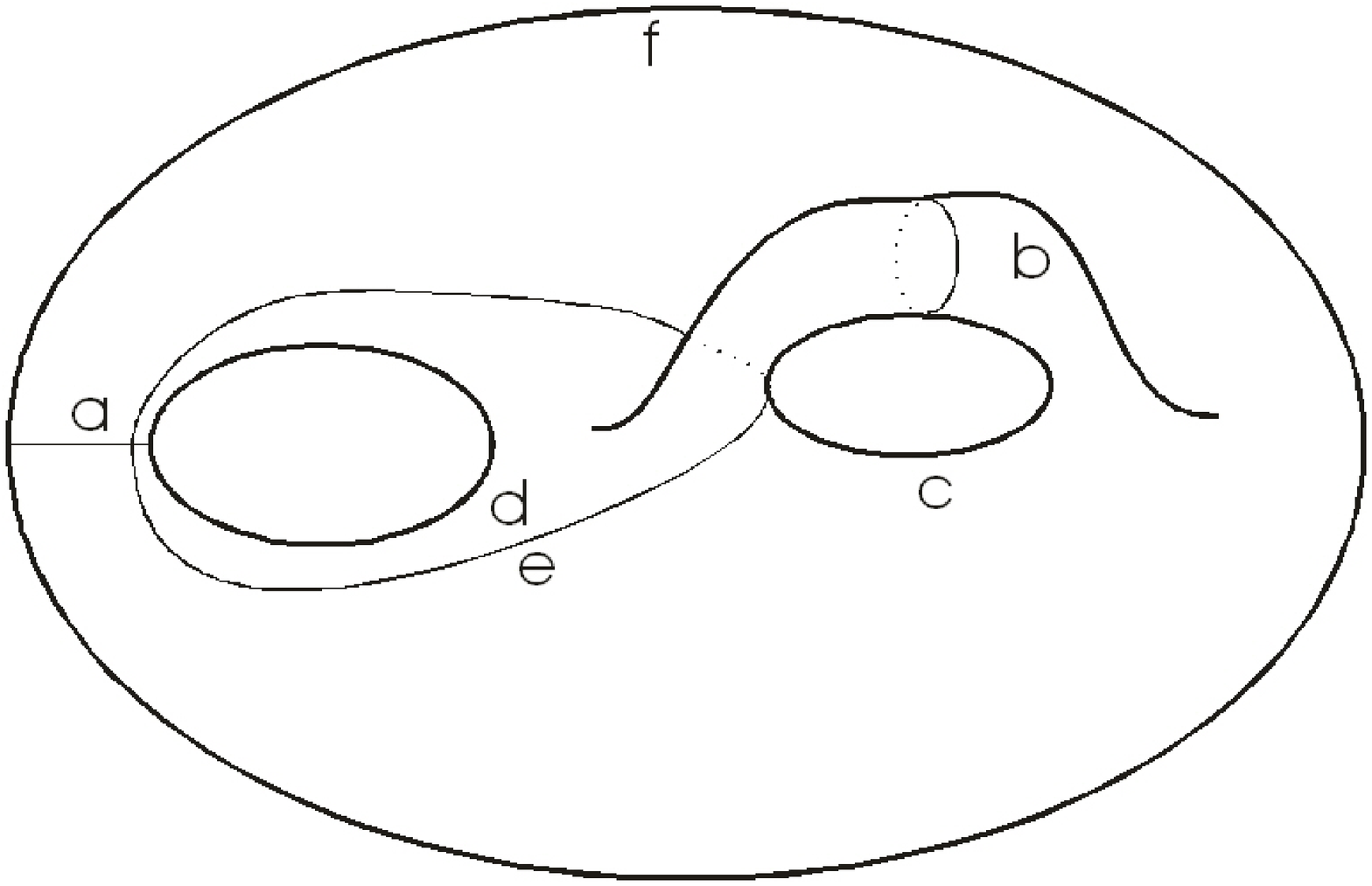}%
\\
\textbf{Figure 3: Curve lengths represent position in moduli space}%
\end{center}
%EndExpansion

\end{center}

Returning to our annular domain metric structure varying with fixed compact
range and boundary conditions from, table 1. As one length becomes longer in
the domain, we care less about the length of its image, as the corresponding
Jacobian entries will be small anyway. As it becomes shorter then we care more
about its length in the image, as we can make big savings on energy (integral
of sum of squares of Jacobian entries) by reducing the value of the
corresponding Jacobian entries when their values are high. So when the domain
becomes a thin cylinder we care about shrinking the circles in the image and
do not care about the length of curves that connect boundary components. Hence
the image collapses down to line segments along the axis of rotation union
discs at the boundary components.

Similarly when the domain becomes a thin ribbon, almost a circle, the distance
between boundary components in the domain is very short so we care about
making the distance between boundary components in the image short. The
corresponding Jacobian entries will be large, because the boundary components
in the range are a fixed distance apart. Conversely we care less about the
length of the meridian circles in the image. Hence in the limit the images
become ruled surfaces in this example with geodesics in $\mathbb{R}^{3}$
connecting the boundaries. The fact that these entries must go to infinity
when the domain distances become small makes the energy go to infinity at this
type of limit of moduli space.

To summarize, in limiting behaviors, proportionately longer curves in the
domain get mapped to longer curves in the image and proportionately shorter
curves in the domain get mapped to shorter curves in the image that can
collapse to a point if they bound a disc in the range.

\subsection{Defining a common domain for taking limits of maps}

We need topologies for the limits of domains, maps and images as
$a\rightarrow0$ and $a\rightarrow\infty.$ For the domain limits as shown in
table 1, we cannot use the Hausdorff set topology inherited from the metric of
the intervals $\left[  -a^{2},a^{2}\right]  \mathsf{X}\left[  0,2\pi\right]
\in\mathbb{R}^{2}$ because these sets are unbounded. Instead we can use it if
we take a uniformly compact representative of each conformal class such as
$\left[  -a^{2},a^{2}\right]  \mathsf{X}\left[  0,2\pi\right]  $ for
$a\rightarrow0$ and $\left[  -1,1\right]  \mathsf{X}\left[  0,2\pi
/a^{2}\right]  $\ for $a\rightarrow\infty$.

For technical reasons we need to treat these domains as a fixed annulus
[-1,1]$\mathsf{X}$[0,1] with a varying metric;

$g\left(  a\right)  _{ij}=\left[
\begin{array}
[c]{cc}%
a^{2} & 0\\
0 & 2\pi
\end{array}
\right]  g_{ij}$, for $a\rightarrow0$ and $g\left(  a\right)  _{ij}=\left[
\begin{array}
[c]{cc}%
1 & 0\\
0 & \frac{2\pi}{a^{2}}%
\end{array}
\right]  g_{ij}$, for $a\rightarrow\infty$

Denote the domain limits $d_{0}$ and $d_{\infty}$\ respectively. We can now
define the matrices as deformations $f_{i}$ such that $d_{i}=f_{i}(d_{1})$.
This defines $h_{i}\circ f_{i}$\ the pull-back maps of the harmonic maps to
the fixed domain $d_{1}$, [-1,1],[0,2$\pi$], enabling point-wise limits to be
taken. When bubbling occurs the point-wise limit map's image will be a subset
of the Hausdorff set limit of the images.

\subsection{Domains represent conformal classes}

However, to be precise, we are really interested in the map $H$ from the space
of conformal classes $C_{a}$ in moduli space $D$ of the domain to an image set
$X$ in the range. Here $a$ is the $a$ in table 2 which indicates conformal
class. So $image(h(d_{a}))$ is the image set of a harmonic map $h$ from metric
representative $d_{a}$ of conformal class $C_{a}$. The image and energy of $h$
is an invariant of $C_{a}$.

\begin{definition}
$H:C_{a}\rightarrow X\Leftrightarrow\forall d_{a}\in C_{a},$ $image(h(d_{a}%
))=X$
\end{definition}

In this case the limit would result from combining the above domain and range
limits, enabling us to extend $H$ to $C_{0}$ and $C_{\infty}$, the boundary
points of moduli space where $a=0$ and $a=\infty$ in section 2.1.

\begin{definition}
$\overline{H}:C_{a}\rightarrow X\Leftrightarrow\underset{a\rightarrow i}{\lim
}$ $\left(  image(h(d_{a}))\right)  $ for $i\in\lbrack0,\infty]$
\end{definition}

\subsection{A topology for image limits.}

In this section we outline a new technique [MS1] for taking limits of surfaces
using geometric measure theory on the sphere bundle of $\mathbb{R}^{3}$. This
gives us rectifiable sets as limits of rectifiable sets and captures lower
dimensional sets that arise in the limit such as those resulting from surfaces
collapsing down curves or points. Consider the canonical example of a sequence
of circles in $\mathbb{R}^{2}$ of radii tending to zero. Their lift in the
sphere bundle will tend to a line, see figure 4. This projects down a positive
measure concentrated at a point.

\begin{center}%
%TCIMACRO{\FRAME{dtbpFU}{1.1366in}{1.1424in}{0pt}{\Qcb{\textbf{Figure 4: A
%sequence (below) and its lift (above)}}}{}{Figure}%
%maintain-aspect-ratio TRUE;  display "USEDEF";  valid_file "T";
%width 1.1366in;  height 1.1424in;  depth 0pt;  original-width 1.0693in;
%original-height 1.0751in;  cropleft "-0.075164";  croptop "1.014953";
%cropright "0.924836";  cropbottom "0.014953";
%tempfilename 'I3LT5P0A.wmf';tempfile-properties "XPR";}}}%
%BeginExpansion
\begin{center}
\includegraphics[
trim=-0.080373in 0.016076in 0.080373in -0.016076in,
%natheight=1.075100in,
%natwidth=1.069300in,
%height=1.1424in,
width=1.1366in
]%
{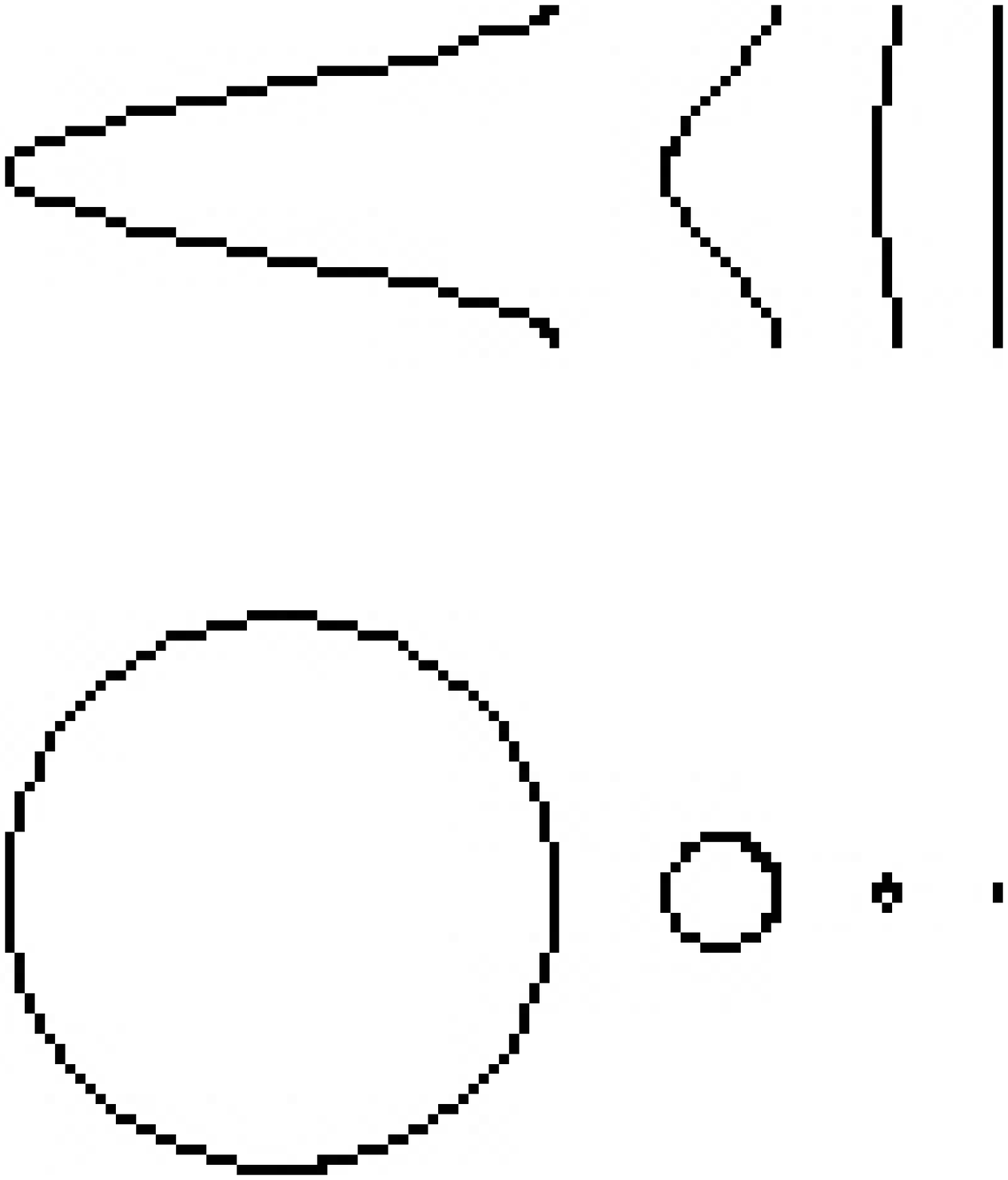}%
\\
\textbf{Figure 4: A sequence (below) and its lift (above)}%
\end{center}
%EndExpansion

\end{center}

\subsection{Image and energy variation}

We recall the qualitative observation on how deforming the domain to different
limits in moduli space will either make energy finite or infinite. For this
section we will consider domains d$_{a}$=[-a,a]\textsf{X}[0,$\frac{2\pi}{a}$].

\subsubsection{\textbf{High energy }$\left(  a\rightarrow0\right)  $.}

We keep domain area constant, so energy depends on the average sum of squares
of entries of the Jacobian. When $a\rightarrow0$, the two boundaries of the
domain annulus become very close. This requires a very high energy for the
image to span the boundary conditions in $\mathbb{R}^{3}. $

\subsubsection{\textbf{Low energy }$\left(  a\rightarrow\infty\right)  $}

Conversely as $a\rightarrow\infty$, the boundaries become far apart in the
domain, enabling energy to be minimized or to stay finite as the Jacobian
entries will mostly be small. However this finite energy limit can exhibit
bubbling, and the image will collapse down to a line segment union one or two discs.

Notice that in the low energy limit in table 1, the energy limits equal twice
the area of the minimal surfaces in the corresponding limits of the images.
This motivates extending the classical energy-area inequality to the limits of
maps. Note also that the inequality can be expressed in terms of $H $ above,
as well as $\overline{H}$, making its statement extendable.

Note that the energy minimizer under domain deformation may or may not be at
the boundary of moduli space, depending on the boundary conditions in the
range. For example if the minimal surface catenoid connecting two boundary
components has lower area that the two discs, then the conformal structure on
the interior of moduli space that maps to the catenoid by a harmonic map will
yield the energy minimizer.

Also note that in the finite energy limit the tubes which collapse down to
straight line segments do not contribute to the energy as they have zero area.
This is of interest as the next section shows that the entries of the Jacobian
that contribute to energy will go to infinity in a thin tube.

This allows a generalization of Douglas' result for existence of certain
minimal surfaces as images of harmonic maps from discs. Douglas takes the
minimum energy map over a disc with fixed conformal structure as boundary
values vary. If instead we allow the domains to have topology and hence non
trivial moduli space, we can find the energy minimizer for maps to multiple
boundary conditions, as conformal structure of domain and boundary values
vary. The thin tubes with zero energy in the limit allow for the image limit
to be unions of minimal surfaces, with separate boundary conditions, connected
by straight line segments.

\subsubsection{\textbf{Energy and image curvature}}

Note that at the `high energy maps' end of moduli space, in our example, the
energy goes to infinity as the image approaches a cylinder or a cone. These
have one principle curvature positive finite and the other zero. This suggests
a relationship between ratios of principle curvatures, and energy of harmonic
maps, as they can both tend to infinity together, giving ruled surface image
limits. Similarly at the low energy end we have minimal surface image limits,
and principle curvature ratios =1. Thus both quantities are minimized together.

This observation, that the ratio of principle curvatures in an image can lead
to infinite energy maps, implies that the inverse problem mentioned in the
introduction may not always have a solution. We can now quote a physical
intuition and result from [MS2]

\subsection{Physical intuition II: Image curvature and energy}

Also we can say that as the image deviates more from being a minimal surface,
more energy is required of the harmonic map. A physical interpretation is to
consider an elastic sheet with low curvature in one direction and high
curvature in the other. To maintain equilibrium, the tension in the sheet in
the low curvature direction must be much greater than in the other direction,
thus contributing more to energy. This effect can be seen in the infinite
energy ruled surface cases, cone and cylinder, in table 1, where the principle
curvature ratios go to infinity.

This can also be extended to the case of a planar sheet with higher in tension
in one direction than in another. If we consider small deformations of the
boundary data, then small curvatures will be introduced into the sheet and
their ratios will reflect the ratios in the tensions in different directions.

\subsection{Improved energy-image curvature inequality}

\begin{theorem}
If h is a degree 1 harmonic map from a smooth compact surface into
$\mathbb{R}^{n}$, and $\rho_{1}$ and $\rho_{2}$ are principle curvatures of
the image of h, then%

\begin{equation}
\underset{}{Energy\geq\underset{}{\underset{image}{%
%TCIMACRO{\diint }%
%BeginExpansion
{\displaystyle\iint}
%EndExpansion
}\left(  \sqrt{\frac{\rho_{1}}{\rho_{2}}}+\sqrt{\frac{\rho_{2}}{\rho_{1}}%
}\right)  d_{i}d_{j}\geq2(area}of}image)
\end{equation}

whenever the integral makes sense on the image taking $\frac{0}{0}=1.$

Furthermore when the pull back to the domain of the directions of principle
curvatures are defined, we can say that equality is achieved on the left hand
side if and only if the pull back of the directions of principle curvatures
are orthogonal in the domain. This occurs for the radially symmetric case.
\end{theorem}

When $\rho_{1}$=$\rho_{2}$ , the right hand side becomes twice the area of the
image. For proof of the inequality see [MS2].

\subsection{Examples of moving bubbles around}

Here are four parameterizations in terms of $s$ and $\theta$, of the annuli
under domain deformations with the metrics inherited from Euclidean space in
the obvious way. The point-wise image limit is shown in Table 3. The bubble in
each case being the discs union the straight line segment minus the image
shown. Each domain is parameterized in terms of a fixed annulus $A=\{(s,\theta
):(s,0)=(s,2\pi),-1\leq s\leq1,0\leq\theta\leq2\pi\}.$ Now first we have the
parameterization of the domain $D_{\varepsilon}=u_{\varepsilon}(A).$ The
harmonic map $h:D_{\varepsilon}\rightarrow\mathbb{R}^{3}$. Note that $u$ is
always uniformly bi-Lipschitz for each $\varepsilon$.

1) D$_{\varepsilon}$ is an identified rectangle in $\mathbb{R}^{2}$ plane,
given by $u_{\varepsilon}:(s,\theta)->(x,y),$ $x=s,$ $y=\varepsilon\theta.$

Here the point-wise limit of $h(u_{\varepsilon}(s,\theta))\rightarrow
(R,\Theta,Z)$ as $\varepsilon\rightarrow0$, the pull backs of the harmonic
maps as $\varepsilon\rightarrow0$ is given by $f(s,\theta)\equiv\left\{
\begin{array}
[c]{c}%
(1,\theta,0)\text{ }for\text{ }s=0\\
(0,\theta,s)\text{ }for\text{ }0<s<1\\
(1,\theta,1)\text{ }for\text{ }s=1
\end{array}
\right\}  $where the mass (two dimensional Hausdorff measure) of the domain is
mapped to the line segment, and the disc interiors bubble.

2) D$_{\varepsilon}$ is the planar annulus. In polar coordinates
$r=\varepsilon+s(1-\varepsilon)$, $z=0$. Here the point-wise limit of the pull
backs of the harmonic maps as $\varepsilon\rightarrow0$ is given by
$f(s,\theta)\equiv$ $\left\{
\begin{array}
[c]{c}%
(1,\theta,0)\text{ }for\text{ }s=0\\
(1,\theta,1)\text{ }for\text{ }s>0
\end{array}
\right\}  $

Here the bubble is one boundary ring, the interiors of the two discs and the
straight line segment.

3) The double cone/hyperboloid with a metric from cylindrical coordinates in
$\mathbb{R}^{3}$ given by $z^{2}=r^{2}+\varepsilon,-1<z<1$. Here the
point-wise limit of the pull backs of the harmonic maps as $\varepsilon
\rightarrow0$ is given by $f(s,\theta)\equiv\left\{
\begin{array}
[c]{c}%
(1,\theta,0)\text{ }for\text{ }s<0.5\\
(0,\theta,0.5)\text{ }for\text{ }s=0.5\\
(1,\theta,1)\text{ }for\text{ }s>0.5
\end{array}
\right\}  $

Here the bubble is the interiors of the two discs and the straight line
segment minus the center point.

4) The spherical annulus. That is an annulus in spherical coordinates with
$\rho=\varepsilon+s(\pi-\varepsilon)$,with metric inherited from
$\mathbb{R}^{3}$. Here the point-wise limit of the pull backs of the harmonic
maps as $\varepsilon\rightarrow0$ is

$f(s,\theta)\equiv\left\{
\begin{array}
[c]{c}%
(1,\theta,0)\text{ }for\text{ }s=0\\
(0,\theta,0.5)\text{ }for\text{ }0<s<1\\
(1,\theta,1)\text{ }for\text{ }s=1
\end{array}
\right\}  $

Here the bubble is the interiors of the two discs and the straight line
segment minus the center point.

To derive these results consider a conformal map from each $D_{\varepsilon}$
to an annulus in the shape of part of the surface of a cylinder. This involves
the identity on $\theta$, and scaling the longitudinal directions along the
surface by $\frac{1}{r}$ to give a conformal map. Then $h$ will map the
longitudinal direction linearly onto $Z$. Equation 2.2 then gives the $R(Z)$
coordinate, with $\Theta=\theta$.

\begin{center}%
%TCIMACRO{\FRAME{dtbpFU}{3.4205in}{3.4205in}{0pt}{\Qcb{Table 3: Examples of
%bubbling}}{}{Figure}{\special{ language "Scientific Word";  type "GRAPHIC";
%display "USEDEF";  valid_file "T";  width 3.4205in;  height 3.4205in;
%depth 0pt;  original-width 5.9867in;  original-height 4.8418in;
%cropleft "0";  croptop "1.0011";  cropright "0.8095";  cropbottom "0";
%tempfilename 'I3LT5P0B.wmf';tempfile-properties "XPR";}}}%
%BeginExpansion
\begin{center}
\includegraphics[
trim=0.000000in 0.000000in 1.140466in -0.005326in,
%natheight=4.841800in,
%natwidth=5.986700in,
%height=3.4205in,
width=3.4205in
]%
{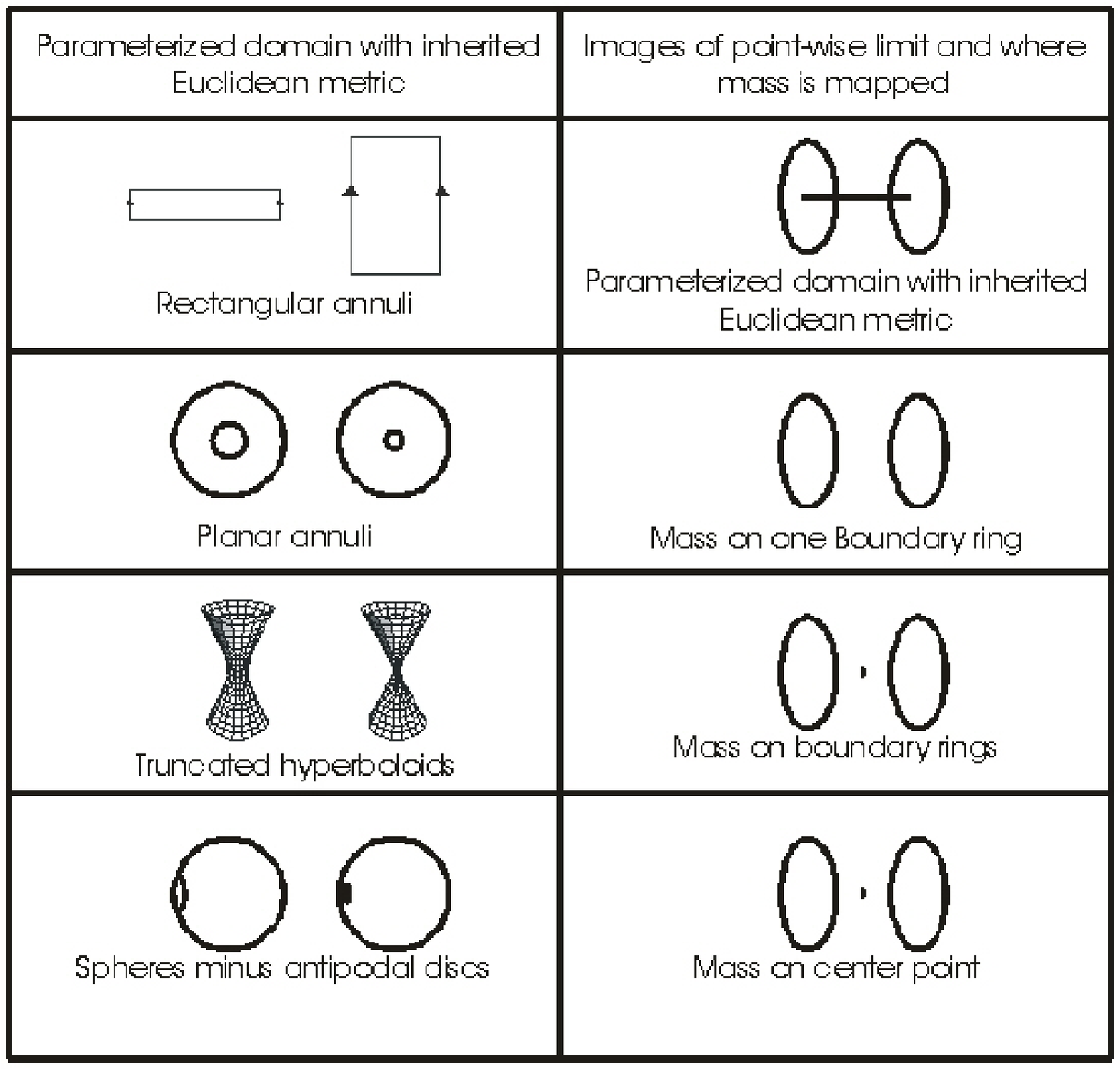}%
\\
Table 3: Examples of bubbling
\end{center}
%EndExpansion

\end{center}

We can move the bubble and discontinuity around, from row to row in table 3,
by changing the metric representatives of the conformal classes of domain. The
next section uses this flexibility to eliminate bubbling by making the domain
metric limit close to the limit image metric, and invoking a compactness
theorem [MS1] applicable to sets with the regularity of images of harmonic
maps in Riemannian manifolds.

\section{Eliminating bubbling in harmonic map sequences.}

Bubbling is the phenomenon of limits of maps where, if we consider the limit
of the graphs of the maps, a region of the image lies over a lower dimensional
region of the domain. In the case of sequences of harmonic maps with varying
domains we have to set up a fixed domain to be able to talk about bubbling. In
section 2.3 we composed the harmonic maps with the unique linear map from the
identified unit square to the domains of the harmonic maps.

In our examples from table 1 bubbling is always associated with a
discontinuity in the point-wise limit. In general this need not happen with
bubbling, but bubbling cannot occur without derivatives becoming unbounded as
relatively smaller parts of the domain are mapped to larger parts of the
image. Consider the limit of functions: $\underset{n\rightarrow\infty}{\lim
}\left(  nxe^{-\left(  nx\right)  ^{2}}\right)  $. The point-wise limit is
zero everywhere, despite there being a bubble.

However we will see in certain cases that if we construct a sequence of domain
metrics that tend to the limit in the same way as the metric structures of the
images, then we do not have bubbling. First we will consider more general
cases before giving an outline of the proof for the radially symmetric case.

\subsection{Physical intuition III: Maps near the identity do not bubble}

Aiming to make the limit of maps like the identity map may minimize energy in
some sense and eliminate bubbling, by aiming to keep derivatives bounded near
to 1.

\subsection{Constructing metric structures on domains}

Allow the domains D$_{n}$ and Images I$_{n}$ to be topologically equivalent;
e.g.: an annulus as in section 2. Then we can define a condition to avoid
bubbling by comparing the metrics on the domains and images. To enable
comparison, we set up pull back metrics from both domains and ranges on a
topologically equivalent manifold $M$ via the maps, $\Psi_{D_{n}}:M\rightarrow
D_{n}$ and $\Psi_{I_{n}}:M\rightarrow I_{n}$ giving $d_{D_{n}}(x,y)=d(\Psi
_{D_{n}}(x),\Psi_{D_{n}}(y)),\forall x,y\in M$\ for the domains and similarly
for the images. Thus the no-bubbling condition can be defined as:

$\underset{n\rightarrow\infty}{\lim}d_{D_{n}}(x,y)=\underset{n\rightarrow
\infty}{\lim}d_{I_{n}}(x,y),\forall x,y\in M$\qquad\qquad(Condition 3.1)

Note that in the dimension collapsing case these limits can be zero. Bubbling
involves unbounded distance ratios between points in the domain and the range.
This is impossible in the above condition. The following lemmas show that the
limits of the graphs will be the identity, in suitable coordinate systems.
Lemmas 3.1 and 3.2 develop this idea, and allow us to state theorem 3.3 which
is developed further in special cases in theorem 3.4.

\begin{lemma}
(special case)

The identity map from a surface to a surface is harmonic, when the surface is
the range.
\end{lemma}

\begin{proof}
This is because energy will equal twice the area, giving equality for the
energy inequality, and achieving the minimum possible.
\end{proof}

\begin{lemma}
If two sequences of metrics, $d_{1,i}$ and $d_{2,i}$, on a surface are close
in the sense that a bi-Lispchitz bijection exists between them with Lipschitz
constants less than $1+\varepsilon,(0<\varepsilon<1/i)$, then such bijections
will have energy that tends to twice the area of the image as $\varepsilon
\rightarrow0$.
\end{lemma}

\begin{proof}
From the bi-Lipschitz condition we can write:

$\frac{(\text{\textit{area of image}})}{(1+\varepsilon)^{2}}\leq($\textit{area
of domain}$)\leq($\textit{area of image}$)(1+\varepsilon)^{2}.$

Now $Energy\leq2($\textit{area of domain}$)(1+\varepsilon)^{2}$ by definition
of energy. These combine to give $Energy\leq2($\textit{area of image}%
$)(1+\varepsilon)^{4}$ . As $\varepsilon\rightarrow0$ we approach the desired
equality, as $Energy$ $\geq$ $2($\textit{area of image}$)$ for any map of surfaces.
\end{proof}

As the image only depends on the conformal class of the domain, we can then
choose the metric representatives, $d_{D_{n}}$\ in the sequence to satisfy
condition 3.1. We can now conclude:

\begin{theorem}
Suppose the hypotheses for lemma 3.2 are satisfied by a sequence of domains
and images of a sequence of harmonic maps which themselves need not be
uniformly bi-Lipschitz. Then we can set up a sequence of uniformly
bi-Lipschitz maps with energy and image tending to the same limits as the
energies and images of the original harmonic map sequence.
\end{theorem}

The proof follows from Lemma 3.2. We will prove a stronger theorem in the
special case of our radially symmetric examples from section 2.

\begin{theorem}
When the hypotheses for lemma 3.3 are satisfied in the radially symmetric case
with annular domains the new sequence in the conclusion of theorem 3.3 can be
chosen to be a sequence of harmonic maps.
\end{theorem}

We construct a set of domains that approaches the limit image, as suggested by
thm. 3.2. See figure 5. Two annuli connected by a thin tube are given
dimensions to fit the boundary data in $\mathbb{R}^{3}.$ As $\varepsilon
\rightarrow0$, the boundary of moduli space is approached and the domain
approaches the image limit.

\begin{center}%
%TCIMACRO{\FRAME{dtbpFU}{2.8119in}{1.7617in}{0pt}{\Qcb{\textbf{Figure 5:
%Domains to eliminate bubbling}}}{}{Figure}%
%maintain-aspect-ratio TRUE;  display "USEDEF";  valid_file "T";
%width 2.8119in;  height 1.7617in;  depth 0pt;  original-width 5.095in;
%original-height 3.1747in;  cropleft "0";  croptop "1";  cropright "1";
%cropbottom "0";  tempfilename 'I3LT5P0C.wmf';tempfile-properties "XPR";}}}%
%BeginExpansion
\begin{center}
\includegraphics[
%natheight=3.174700in,
%natwidth=5.095000in,
%height=1.7617in,
width=2.8119in
]%
{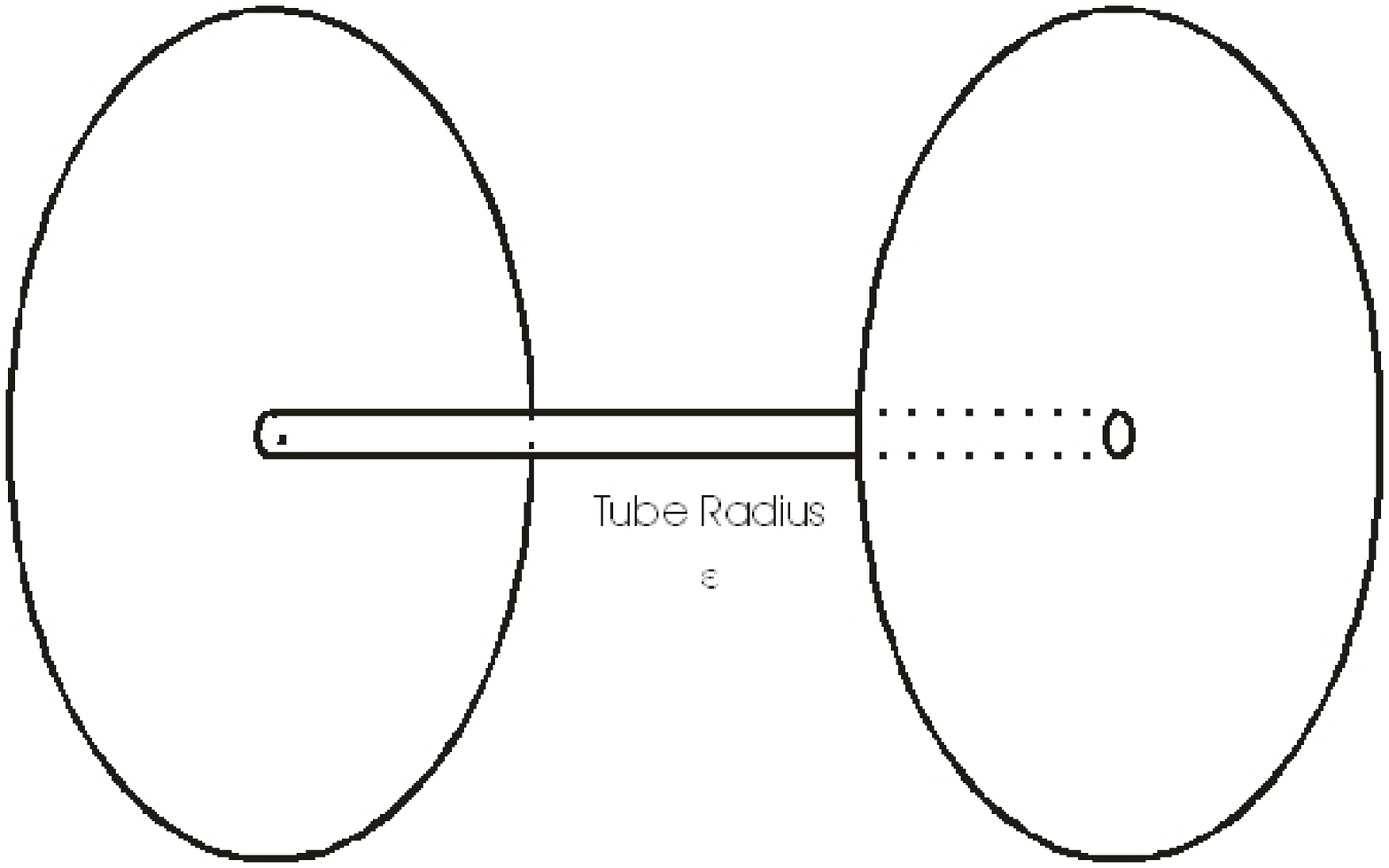}%
\\
\textbf{Figure 5: Domains to eliminate bubbling}%
\end{center}
%EndExpansion

\end{center}

\begin{proof}
Given symmetric boundary data we know that the image will satisfy equation
2.2. As we approach the boundary of moduli space in the domain we can show
that in our approximation to the image, the tube is mapped mainly to the tube
and the punctured discs are mapped mainly to the punctured discs. We do this
by making energy calculations for small thin annuli near the puncture of the
discs in figure 5.

Part 1: Tubes do not fan out into the annuli.

Consider an annulus in the domain of inner radius $r,$ and outer radius $2r$
centered on the axis. Lets say it maps to a radius of at least $k.$The map has
a tangential derivative of at least $\frac{k}{2r}$, and so the energy for this
small annulus is at least $\left(  \frac{k}{2r}\right)  ^{2}3\pi r^{2}%
=\frac{3\pi k^{2}}{4}$. However as $\varepsilon\rightarrow0$, we need an
infinite number of such annuli of radii $\frac{r}{2^{n}},$ so the energy
becomes infinite, unless there is no such $k$%
%TCIMACRO{\TEXTsymbol{>}}%
%BeginExpansion
$>$%
%EndExpansion
0 that the tubes fan out to.

Part 2: The annuli do not stretch into the tubes.

We calculate the image of the neck of the tubes. First we need to conformally
map a domain such as that in figure 5 to an identified rectangle annulus as in
table 1. Let this annulus be a cylindrical surface of circumference 2$\pi,$
and let the domain, without loss of generality, have length 1 and boundary
radii 1. The conformal map will scale the length along the cylinder by the
reciprocal of the radius of the domain.

Each planar annulus in figure 5 will be mapped to a length of $\left\vert
log(\varepsilon)\right\vert .$ This is explained in the final paragraph of
section 2.9. The tube part in figure 5 will be mapped to a length of $\frac
{1}{\varepsilon}.$ We know from section 2.1 that $Z$ (varies from -1 to 1) is
a linear function of $x$ (varies from $-\left\vert log(\varepsilon)\right\vert
+\frac{1}{2\varepsilon}$ to $\left\vert log(\varepsilon)\right\vert +\frac
{1}{2\varepsilon}$), so we can place the image of one of the necks at
$1-\left(  \frac{\left\vert log(\varepsilon)\right\vert }{\left\vert
log(\varepsilon)\right\vert +\frac{1}{\varepsilon}}\right)  .$ As
$\varepsilon\rightarrow0,$ this becomes 1. Therefore the images of the necks
remain at the ends of the tube as required. We have shown that the sequence of
harmonic maps tend to the identity and so do its derivatives. It satisfies
condition 3.1, therefore it does not bubble.
\end{proof}

\begin{conjecture}
Whenever the domain deformation and induced image variation tend to the same
boundary point of moduli space we can eliminate bubbling by choice of domain
representatives of conformal classes, and if needed by choice of harmonic map.
[MS1] gives a topology in which images of harmonic maps converge to give the
desired limit image with dimension collapsing. This limit image can then be
used to construct a set of domains which will lead to a sequence of harmonic
maps which converge point-wise to the limit image.

\begin{theorem}
For the construction in thm 4.1 conjecture 1 is true

\begin{proof}
As proof for thm 3.4
\end{proof}
\end{theorem}
\end{conjecture}

\section{Generalizing to more boundary components}

\subsection{Minimal surfaces union straight line segments}

The following theorem was conjectured by Robert Hardt.

\begin{theorem}
For connected planar (genus zero) multiply connected domains with more than 2
boundary components, there will be a sequence of harmonic maps from domains
approaching a boundary point in moduli space whose energies approach the
infimum of all energies of harmonic maps from that topological type of surface
over all conformal structures with range $\mathbb{R}^{n}$ and a union of
suitably embedded $\mathbb{S}^{1}$s as Douglas type boundary conditions.
Furthermore the limit of images under a suitable topology that captures
dimension collapsing such as in [MS1] will be a union of minimal surfaces
interconnected by an embedded complex of straight line segments.
\end{theorem}

We construct a domain consisting of $n$ planar annuli each connected by a thin
tube of length 1 and diameter $\varepsilon.$ These tubes meet in a central
joint, so that the whole domain has the topological type of an $\mathbb{S}%
^{2}$ with n discs removed. The sequence of domains is given by $\varepsilon
\rightarrow0.$ Douglas boundary conditions are used for technical reasons.

\begin{theorem}
We know that every section of two tubes connecting any pair of boundary
components will pinch down to an arbitrarily small diameter at least one point.
\end{theorem}

\begin{proof}
This is because if it does not, the energy on that section of tubes will grow
to the order of $\frac{1}{\varepsilon}$ as $\varepsilon\rightarrow0$.
\end{proof}

\begin{proposition}
This means at least n-1 tubes contain a pinch, even if they all pinch at the joint.
\end{proposition}

\begin{proposition}
The images of harmonic maps as $\varepsilon\rightarrow0$ contain annuli
contained in the convex hull of each boundary component and a circle diameter
$\varepsilon$ embedded in $\mathbb{R}^{n}$.
\end{proposition}

\begin{proof}
The convex hull property of images of harmonic maps applies between the neck
pinch and the boundary component.
\end{proof}

\begin{proposition}
As $\varepsilon\rightarrow0$ the image area measure will tend to the sum of
the areas of minimal surfaces each with one of the prescribed boundary components
\end{proposition}

\begin{proof}
For a given conformal class with $\varepsilon$ small we can show that the area
of the annular part of the image is within a fixed range of the area $A$ of
the minimal surface bounded by the boundary component. We do this by
constructing a candidate image with $energy$\ $\leq2A+e(\varepsilon)$, where
$e$ is arbitrarily small. We know by the area-energy inequality that the area
of the actual image $A\prime$ is within $e/2$ of $A$. The candidate images are
the minimal surfaces each with a disc removed (from Douglas' result [D]) and a
tube of diameter $\varepsilon$ attached leading to a common ball of diameter
$2\varepsilon$ where the tubes join with the required topology.
\end{proof}

\begin{proof}
(of \textbf{theorem 4.1}) A deviation from the minimal surface can only be
achieved by thin immersed tubes. The topology of the domain ensures that such
thin tubes will be retractable, in an appropriate sense for immersed tubes, to
the surface, apart from the ones required by domain topology. The retraction
in each case can be done so as to decrease energy of the map. The tubes will
be straight line segments by the convex hull property of images of harmonic maps.
\end{proof}

\begin{conjecture}
It may even be possible that there is only one boundary point of a suitably
compactified moduli space which has a finite energy limit of map energies for
a given set of minimal surfaces spanning boundary components.
\end{conjecture}

We need to specify the set of minimal surfaces, to apply theorems, as there
may be more than one. In particular consider the annulus example. In addition
to the two discs there can be one or even two catenoids with the same boundary
conditions. In the annular examples these catenoids do not occur as images at
the boundary of moduli space, but in the higher boundary component case they can.

We can observe that in these cases, of higher dimensional moduli space we can
consider not only the boundary point of moduli space approached by a sequence
of domains but also the direction of approach of the sequence within higher
dimensional moduli space. We should also recall that the choice of metric
representative of conformal class can be used to eliminate bubbling, by use of
domains with both discs and thin tubes.

\subsection{Direction of approach to the boundary of moduli space}

\begin{center}%
%TCIMACRO{\FRAME{itbpF}{1.1631in}{1.2038in}{0in}{}{}{Figure}%
%maintain-aspect-ratio TRUE;  display "USEDEF";  valid_file "T";
%width 1.1631in;  height 1.2038in;  depth 0in;  original-width 1.2611in;
%original-height 1.3059in;  cropleft "0";  croptop "1";  cropright "1";
%cropbottom "0";  tempfilename 'I3LT5O00.wmf';tempfile-properties "XPR";}}}%
%BeginExpansion
{\includegraphics[
%natheight=1.305900in,
%natwidth=1.261100in,
%height=1.2038in,
width=1.1631in
]%
{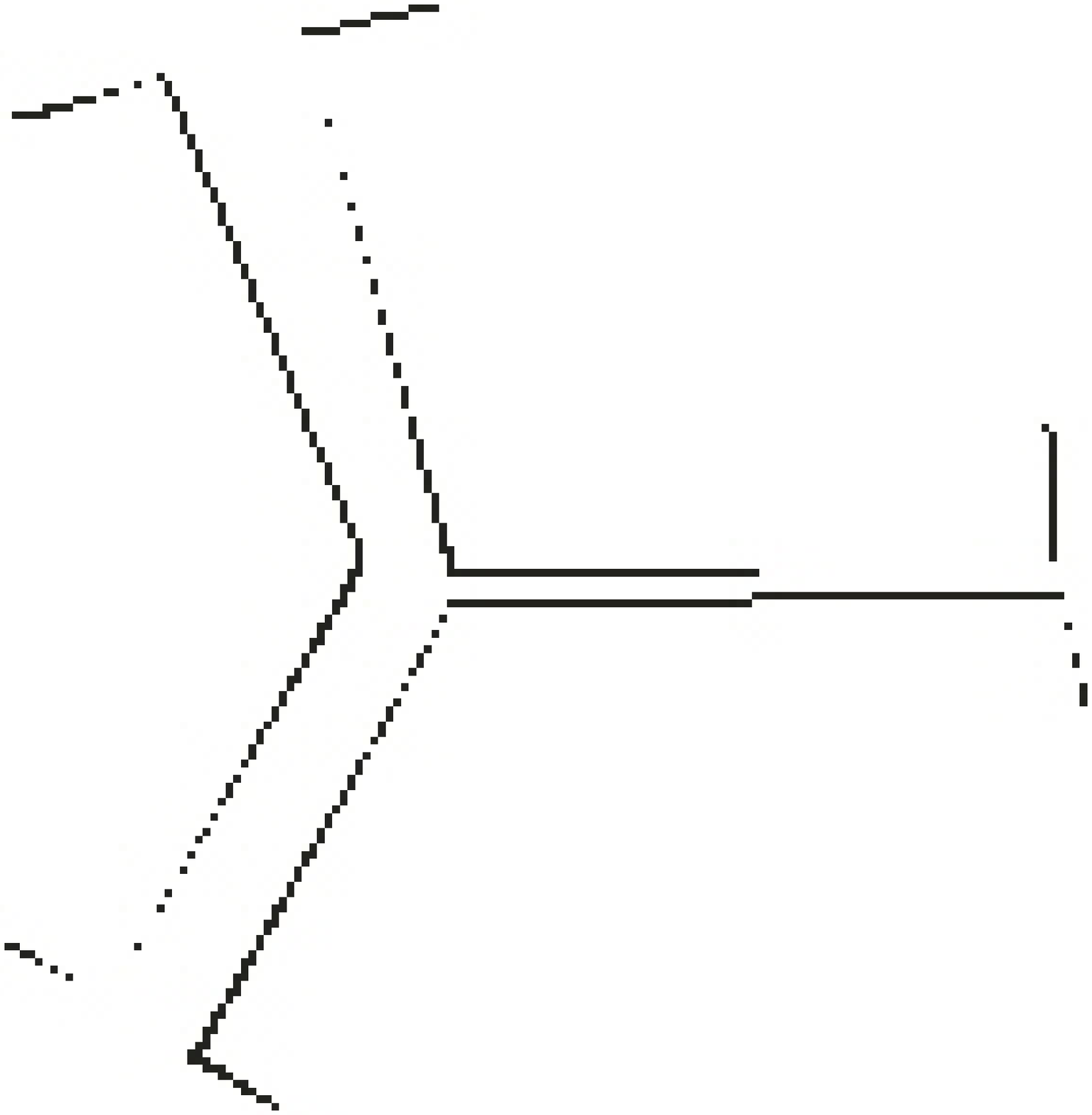}%
}%
%EndExpansion%
%TCIMACRO{\FRAME{itbpFU}{1.1745in}{1.2095in}{0in}{\Qcb{}}{}{Figure}%
%maintain-aspect-ratio TRUE;  display "USEDEF";  valid_file "T";
%width 1.1745in;  height 1.2095in;  depth 0in;  original-width 1.1465in;
%original-height 1.1822in;  cropleft "0";  croptop "1";  cropright "1";
%cropbottom "0";  tempfilename 'I3LT5O01.wmf';tempfile-properties "XPR";}}}%
%BeginExpansion
{\includegraphics[
%natheight=1.182200in,
%natwidth=1.146500in,
%height=1.2095in,
width=1.1745in
]%
{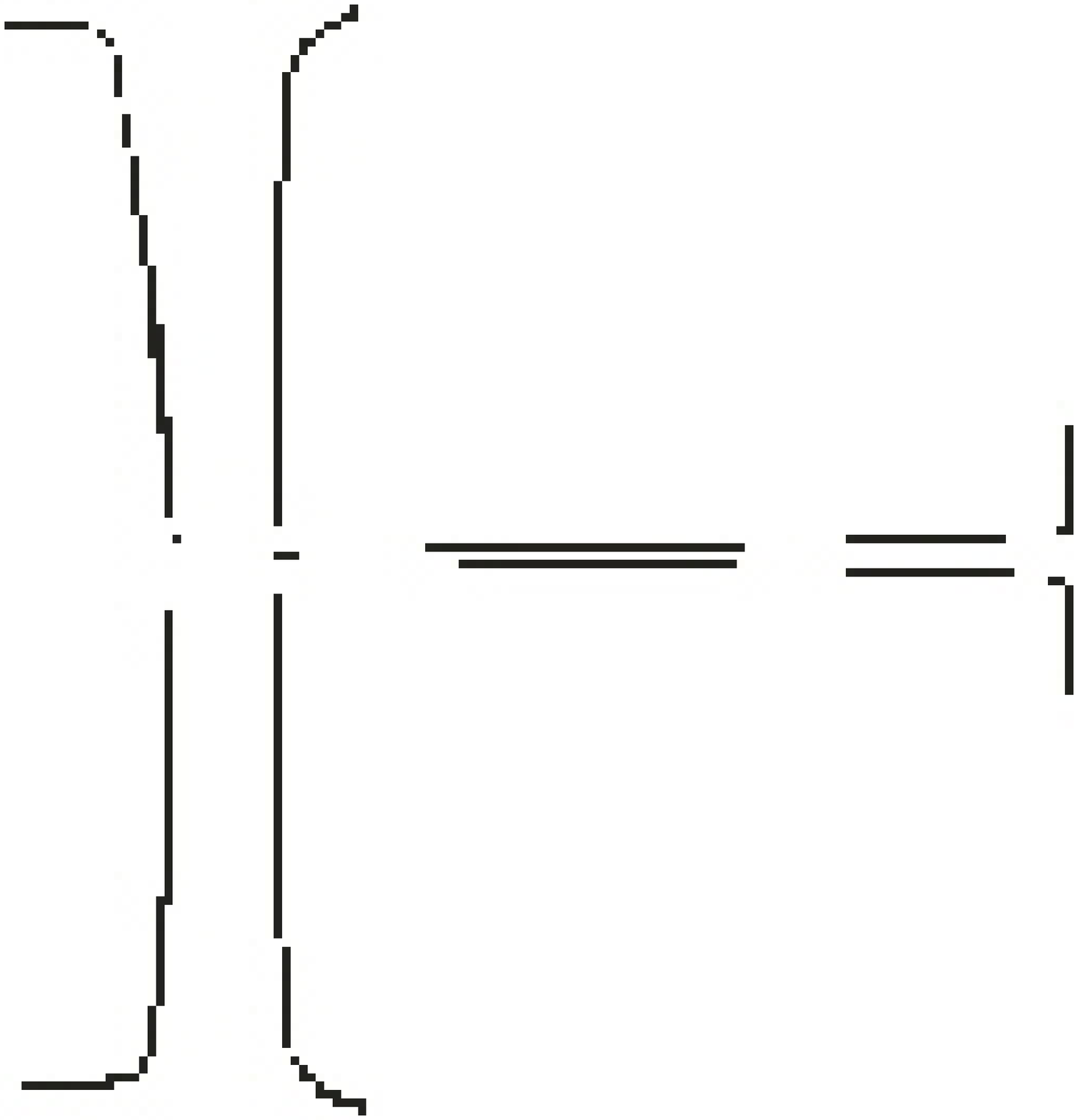}%
}%
%EndExpansion

\textbf{Figure 6 Domain and Image in cross section}
\end{center}

\begin{conjecture}
The surface parts of the limit image will depend only upon the boundary point
of moduli space approached, but the complex of straight line segments can
depend also upon the direction of approach of a sequence, or subsequence with
a convergent image under a topology such as in [MS1].
\end{conjecture}

Consider the example of three boundary components. The moduli space of the
domains will be more than one dimensional. This means that boundary points of
moduli space can be approached by a sequence in different directions on a
surface, so we can investigate if the limit of images of harmonic maps depends
on just the boundary point or the way it is approached within moduli space.

We give a conjectured example in figure 6 that show the direction of approach
to the boundary of moduli space of a sequence of domains making a difference.
We will try to support this with a physical intuition and also a parameterization.

\subsection{\textbf{Physical intuition IV: Force balancing in line segment
skeletons.}}

Three line segments meeting at a point where their forces balance will form
some kind of Y. If all the forces are equal then all the angles will be equal.
If one of the forces is significantly less then the angle in the other two
will lie in $[\pi/3,\pi/2)$.

This models what happens at the finite energy boundary of moduli space. If one
of the domain tubes in figure 6 shrinks in radius more quickly than the rest,
with length constant, then it will act like the line segment with less force.
The limit angle in $[\pi/3,\pi/2)$ will depend upon these ratios, when they
are constant, which are associated with a different path in moduli space. In
this case the limit can be expressed as a stationary varifold. In fact we can
even obtain a T-junction as a limit using a diagonalization argument having
the ratios change. Note that a T-junction cannot be achieved by a three line
segments balancing non-zero forces, i.e. they cannot be represented as a
stationary varifold.

\subsection{\textbf{Parameterized example sequence of domains in moduli
space}}

We suggest the case of the T-junction as a limit of images of harmonic maps
could perhaps be achieved by a sequence of domains that must approach the
boundary of moduli space tangentially. Let $x,y,z$ be the respective radii of
the three tubes of unit length that are all connected at one end in a small
punctured sphere in a Y-junction. We know the moduli space is three
dimensional as the number parameters it has is $3n-6$ where $n$ is the
connectivity or number of boundary components for $n>2$ [Ah, section 5.1]. The
compactified moduli space, is modelled by the positive octant. The origin is
the boundary point to be approached. If we approach on the plane $x=y$ we can
approach along the path $(t,t,kt)$ where $k>0$ and $t$ goes to zero. This will
give a Y-singularity where the angles depend on $k$. If we allow $k$ to vary,
we can let $k\rightarrow0$ as $t\rightarrow0,$ giving us a T-singularity or if
$k$ $\rightarrow\infty$ we simply get two line segments connecting the three discs.

\textbf{Acknowledgement}s

Thanks to Bob Hardt who advised my thesis work, part of which is reported in
this paper, and also to Mike Wolf for his enthusiastic help. Thanks to Bob
Gulliver who helped with supplementary work and guidance for this paper. I
also thank Thierry de Pauw, Chun-chi Lin and Qinglan Xia for their discussions
on geometric measure theory over the years.

\end{document}